\documentclass[11pt]{amsart}
\usepackage{amsmath,amssymb,amscd,amsfonts,verbatim}
\usepackage{paralist}
\usepackage[mathscr]{eucal}

\newtheorem{thm}{Theorem}[section]
\newtheorem{lem}[thm]{Lemma}

\newtheorem{prop}[thm]{Proposition}
\newtheorem{cor}[thm]{Corollary}
\newtheorem{defn}[thm]{Definition}
\newtheorem{assu}[thm]{Assumption}
\newtheorem{assu-nota}[thm]{Assumption--Notation}
\theoremstyle{remark}

\newtheorem{remark}{Remark}

\newcommand{\C}{\mathbb C}

\newcommand{\pp}{\mathbb P}

\DeclareMathOperator{\Alb}{Alb}

\DeclareMathOperator{\Proj}{Proj}
\DeclareMathOperator{\Sym}{Sym}
\DeclareMathOperator{\fdeg}{fdeg}
\newcommand{\epsi}{\varepsilon}

\newcommand{\al}{\alpha}
\newcommand{\be}{\beta}
\newcommand{\ga}{\gamma}

\newcommand{\Ga}{\Gamma}
\newcommand{\De}{\Delta}
\newcommand{\Si}{\Sigma}
\newcommand{\si}{\sigma}

\newcommand{\fie}{\varphi}

\newcommand{\OO}{\mathcal{O}}

\numberwithin{equation}{section}

\title[Severi type inequalities]{Severi type inequalities for irregular surfaces with ample canonical class}
\author{Margarida Mendes Lopes and Rita Pardini}

\thanks{{\it Mathematics Subject Classification (2000)}: 14J29. \\
The first author is a member of the Center for Mathematical
Analysis, Geometry and Dynamical Systems  and the second author is a member of G.N.S.A.G.A.-I.N.d.A.M.  This research was partially supported by the italian  project ``Geometria sulle
variet\`a algebriche e loro spazi di moduli'' (PRIN COFIN 2006) and by FCT (Portugal) through program POCTI/FEDER}
\begin{document}
\begin{abstract} Let $S$ be a smooth minimal complex projective  surface of maximal Albanese dimension. Under the assumption that the canonical class of $S$ is ample and $q(S):=h^0(\Omega^1_S)\ge 5$, we show:
$$K^2_S\ge 4\chi(S)+\frac{10}{3}q(S)-8,$$
thus improving the well known Severi inequality $K^2_S\ge 4\chi(S)$.

We also give stronger inequalities under extra assumptions on the Albanese map or on the canonical map of $S$.

\end{abstract}
\maketitle

\section{Introduction}
Let $S$ be a smooth minimal complex projective  surface of maximal Albanese dimension. The Severi inequality, namely the  relation
$$K^2_S\ge 4\chi(S),$$
has been proven in full generality in  \cite{severi} by means of a limiting argument. A completely different proof,  of a more geometrical nature,  had been  previously given by Manetti in \cite{manetti} under the additional assumption that the canonical class $K_S$ be ample. In \cite{manetti}, the author also  conjectures that the stronger inequality
$$K^2_S\ge 4\chi(S)+4q(S)-12$$
hold for $q(S)\ge 4$. 

Our main result, obtained by arguments of the same type  as in  \cite{manetti}, is a step towards proving this conjecture:
\begin{thm} \label{main}Let $S$ be a smooth surface of maximal Albanese dimension,  with $K_S$ ample and irregularity  $q\ge 5$. Then:
$$K^2_S\ge 4\chi(S)+\frac{10}{3}q-8.$$
\end{thm}

Under extra assumptions on the behaviour of the Albanese map or of the canonical map of $S$, the methods used in proving Theorem \ref{main} yield  better inequalities (notice that in Theorem \ref{main3} it is enough to assume $S$ minimal):

\begin{thm}\label{main3}
 Let $S$ be a smooth minimal surface of maximal Albanese dimension  and irregularity  $q\ge 5$. 
 \begin{enumerate}
\item If  the canonical map  of $S$ is not birational, then 
$$K^2_S\ge 4\chi(S)+4q-17;$$
\item if the canonical map has degree 2 and $S$ has no irrational pencil, then 
$$K^2_S\ge 6\chi(S)+6q-20. $$
\end{enumerate}
\end{thm}

\begin{thm} \label{main2}Let $S$ be a smooth surface of maximal Albanese dimension  and irregularity  $q\ge 5$ with $K_S$ ample. 
\begin{enumerate}
\item If the Albanese map $a\colon S\to A:=\Alb(S)$ is not birational onto its image, then $$K^2_S\ge 4\chi(S)+4q-13;$$
\item  if $S$ has no irrational pencil and  the Albanese map $a\colon S\to A$ is unramified in codimension 1, then 
$$K^2_S\ge 6\chi(S)+2q-8.$$

\end{enumerate}
\end{thm}

\smallskip
The proofs of Theorems \ref{main}, \ref{main3}  and  \ref{main2} are given in \S \ref{proofs}. We remark that Theorem \ref{main2}, (ii) improves the inequalities given in \cite{ca}  and \cite{konnoeven}.

\smallskip
The main tool of this work is a careful analysis, done  in \S \ref{systems}, 
 of two subsystems   of $|K_S|$:
 
1) the subsystem $\pp(\Lambda)$  generated by the zero divisors of the $2$-forms $\al\wedge \be$, for  $\al, \be\in H^0(\Omega^1_S)$,
 
 2) the subsystem $\pp(\Lambda_{\eta})$, generated by the zero divisors of the $2$-forms $\eta\wedge \be$, for  $\eta\in H^0(\Omega^1_S)$ fixed and general and $\be\in H^0(\Omega^1_S)$.
 \smallskip
 
 It is apparent from the proofs given in \S \ref{systems} that the  inequalities above  are by no means sharp. In many  specific situations, for instance for large values of $q$,  the same proofs give sharper inequalities.  However,  in order to  have  unified statements and avoid unnecessarily long proofs, we have chosen  to present our results in the above form.

We also wish to point out that,  if $S$ has an irrational pencil,  then,  under the weaker  assumption that $S$ be minimal,   the slope inequality gives bounds  that  are better than those in Theorem \ref{main} (see 
 \S \ref{PENCIL}).

 \bigskip

{\bf Notation and conventions.}  We work over the complex numbers. All varieties are projective algebraic. A variety has {\em maximal Albanese dimension} if the Albanese map is generically finite onto its image.

As we focus on surfaces, we recall the definition of the numerical invariants of a smooth projective surface $S$: $q(S):=h^0(\Omega^1_S)$ is the {\em irregularity}, $p_g(S):=h^0(K_S)$ is the {\em geometric genus} and $\chi(S):=p_g(S)-q(S)+1$ is the {\em holomorphic Euler characteristic}.  An {\em irrational  pencil} of  a surface $S$ is a fibration  $f\colon S\to B$, where $B$ is a smooth curve of genus $b\ge 2$.  The integer $b$ is called the {\em genus of the pencil}.

Following the most common convention, when dealing with a vector space $V$ we denote by $\pp(V)$ the set of lines of $V$,  while when dealing with a vector bundle ${\mathcal E}$ on a variety $X$ we denote by $\pp({\mathcal E})$ the space $\Proj(\Sym {\mathcal E})$. Hence the fibre of $\pp({\mathcal E})$ over a point $x\in X$ is $\pp({\mathcal E}^{\vee}\otimes \C(x))$. 
\section{Auxiliary results}

\subsection{Pencils on surfaces}\label{PENCIL}

 Recall that by the Castelnuovo--De Franchis Theorem a surface  $S$ has an irrational pencil (i.e. a fibration onto a curve of genus $\ge 2$) if and only if  there are two linearly independent global  $1$-forms $\alpha$ and $\beta$ of $S$ whose wedge product is identically zero.

The following is  a  variation of  the well known Ca\-stel\-nuovo--De Franchis Lemma.
\begin{lem}\label{involution} Let $S$ be a smooth  surface with irregularity  $q\ge 3$, let $\iota$ be an involution of $S$ and let $T$ be a desingularization of $S/i$. 
If $p_g(T)=0$, then $S$ has an irrational pencil.
\end{lem}
\begin{proof} Up to blowing up the isolated fixed points of $\iota$, we may assume that the quotient surface $S/\iota$ is smooth, so that we may take $T=S/\iota$. Since $p_g(T)=0$, the surface $T$ does not have maximal Albanese dimension.  Let $p$ be the irregularity of $S/\iota$. If $p\ge 2$, then the Albanese pencil of $S/\iota$ pulls back to an irrational pencil of $S$. If $p\le 1$, then $q-p\ge2$, hence there exist two linearly independent $1$-forms $\alpha$ and $\beta$ on $S$  such that $i^*\al=-\al$ and $i^*\beta=-\beta$. The $2$-form $\alpha\wedge\beta$ is invariant with respect to $\iota$ and therefore it descends to a global $2$-form of  $T$, which is zero, since $p_g(T)=0$.  So  the form  $\alpha\wedge\beta$ vanishes identically on $S$ and $S$ has an irrational pencil.
\end{proof}

The following result is essentially proven in \cite{xiaopencil}. 
\begin{prop}\label{xiaolines} Let $S$ be a smooth surface of maximal Albanese dimension with irregularity $q$ and let $f\colon S\to B$ be a pencil of genus $b\ge0$ with general fibre $F$. For $\eta\in H^0(\Omega^1_S)$ denote by $\Lambda_{\eta}\subseteq H^0(K_S)$ the image of the map $H^0(\Omega^1_S)\to H^0(K_S)$ defined by $\al\mapsto \eta\wedge \al$.  Let $r\colon H^0(K_S)\to H^0(K_F)$ be the restriction map. 
If  $\eta\in H^0(\Omega^1_S)$ is general, then 
$$\dim r(\Lambda_{\eta})\ge q-b-1.$$
\end{prop}
\begin{proof} Let $U\subseteq H^0(\Omega^1_S)$ be the kernel of the natural map $H^0(\Omega^1_S)\to H^0(K_F)$. 
Consider a subspace $V\subseteq H^0(\Omega^1_S)$ such that $\eta\in V$ and $H^0(\Omega^1_S)=U\oplus V$. 
Then the argument in the lemma on page 599 of \cite{xiaopencil} shows that $r(\eta\wedge V)\subseteq r(\Lambda_{\eta})$ has dimension $\ge q-b-1$.

\end{proof}
\bigskip

As stated in the introduction, in the presence of an irrational pencil,   the slope inequality (\cite{xiaoslope}, cf. also \cite{ch}, \cite{stoppino})   gives  bounds  that  are better than those in Theorem \ref{main}:

\begin{prop} \label{mainpencil} Let $S$ be a smooth minimal surface of  general type of maximal Albanese dimension with  irregularity $q\ge 3$.
Assume that $S$ has an irrational   pencil $f\colon S\to B$, where $B$ is a curve  of genus $b\ge 2$. Then:
\begin{enumerate}
\item $K^2_S\ge 4\chi(S)+4(q-3)$;
\item if  $S$ is not the product of a curve of genus $2$ and a curve of genus $q-2$,  then $K^2_S\ge 4\chi(S)+ 4(q-2)$.
\end{enumerate}
\end{prop}
\begin{proof} 
Assume that $S$ is the product of a curve of genus $b\ge2$ and a curve of genus $q-b\ge 2$. In this case one has:
\[K^2_S=8\chi=4\chi(S)+4(b-1)(q-b-1),\]
and  statements (i) and (ii) are easy to check.

Assume now that $S$ is not a product of curves and denote by $g$ the genus of the general fibre of $f$.
 By the Lemme on page 345 of \cite{appendix}, we have $b+g\ge q+1$.  By \cite[Thm.1]{xiaoslope} we have:
\begin{equation}\label{slope}
K^2_S\ge 4\chi(S)+4(g-1)(b-1)
\end{equation}
and statement (ii) follows easily. 
\end{proof}

\subsection{Linear systems on curves}

The following is a well known refinement of the classical Clifford lemma (see for instance \cite[Lem. 5.1]{beauville} and \cite{laquila}).
\begin{lem}[Clifford+]\label{clifford} Let $C$ be a smooth projective curve of genus $g$ and let $|\De|$ be a linear  system of $C$ such that $\deg \De\le g-1$. If the map given by $|\De|$ is birational, then:
$$\deg \De\ge 3h^0(D)-4,$$

\noindent with equality holding only if $2\De\equiv K_C$.
\end{lem}

Another result that  we need in our proofs is Castelnuovo's bound (cf \cite[\S 1]{harris}).

\begin{lem}[Castelnuovo's bound]\label{castelnuovo}
Let $C$ be a smooth  curve   of genus $g(C)$ and let $f\colon C\to\pp^r$ be a map that is birational onto its image. Denote by $d$ the degree of $f(C)$ and write $d-1=m(r-1)+\epsi$, where $m\ge 0$ and $0\le \epsi\le r-2$ are integers. Then:
$$g(C)\le{ m \choose 2}(r-1)+m\epsi.$$
\end{lem}
\subsection{Surfaces in projective space}
We recall the following result from \cite{beauville} (Lem. 1.4 and Rem.1.5):
\begin{prop}\label{ruled}
  Let $\Si\subset \pp^r$ be an irreducible non degenerate surface of degree $d$. If $d<2r-2$, then $\Si$ is ruled.
  \end{prop}

A ``quantitative'' version of the result above is the following (cf. \cite{milesdeg}, Cor. 1.1):
\begin{prop}
\label{Mileslines}
 Let $\Si\subset \pp^r$ be an irreducible non degenerate surface of degree $d$. If $3d<4r-9$, then $\Si$ is ruled by lines.
\end{prop}

\section{Maps, linear systems and inequalities}\label{systems}

Throughout all the section we make the following 
\begin{assu}
$S$ is  a smooth  surface of maximal Albanese dimension and  irregularity $q\ge 3$.
\end{assu}
Set $A:=\Alb(S)$, let $a\colon S\to A$ be the Albanese map and denote by $\Si$ the image of $a$. We  introduce and study some maps associated to the canonical map  and the Albanese  map of $S$  and the corresponding linear systems.

\subsection{The system $|D|$ and the map $\fie_D$}
We denote by $T$ the tangent space to $A$  at the origin and we set $\pp^{q-1}:=\pp(T)$. 
Let $G$ denote the Grassmannian of lines in $\pp^{q-1}$ and consider  the Gauss map $\ga \colon \Si\to G$, which   maps a smooth point $x\in \Si$ to the line in $\pp^{q-1}$ obtained by  translating  in the origin the tangent space to $\Si$ at $x$  and projectivizing it. 
 \begin{lem} \label{gauss}
If the image of the Gauss map $\ga\colon \Si\to G$ is not a surface, then  $S$ has an irrational pencil of genus $q-1$.
\end{lem}
\begin{proof} Since $\Si$ generates $A$,  by \cite[ch. VIII, Prop. 4.1]{deb2} either $\ga(\Si)$ is a surface or $\ga(\Si)$ is a curve and there is a $1$-dimensional abelian subvariety $K$ of $A$ such that $\Si +K=\Si$.

Assume that we are in the latter  case. The image of $\Si$ in the abelian variety $A/K$ is a curve $B$ that generates $A/K$, hence $B$ has geometric  genus at least $q-1$. On the other hand, $S$ does not admit a map onto a curve of genus $\ge q$ by assumption, hence $B$ has genus exactly $q-1$ and it is smooth. Composing the Albanese map of $S$ with the projection $A\to A/K$ one gets a map $f\colon S\to B$. The fibres of $f$ are irreducible, since otherwise the Stein factorization of $f$ would give a dominant map $S\to B'$, with $B'$ a curve of genus $>q-1$.
\end{proof}

We denote by $\Lambda$ the image of the natural map $v\colon \bigwedge^2H^0(\Omega^1_S)\to H^0(K_S)$ and  by $\bar{v}\colon \pp(\bigwedge^2 H^0(\Omega^1_S))\to \pp(H^0(K_S))$ the corresponding  map of projective spaces. We write  $\pp(\Lambda):=F+|D|$, where $|D|$ is  the moving part. The map $\fie_D$ given by $|D|$ is the composition $\ga\circ a$, followed by the Pl\"ucker embedding of $G$ into $\pp(\wedge^2T)$. Notice that, because of the natural inclusion $|D|\subseteq |K_S|$, the map $\fie_D$ is obtained from  the canonical map of $S$ by  composing it  with a projection.

\begin{lem}\label{fieD} 

\begin{enumerate}
\item If  $S$ has no irrational pencil of genus $q-1$, then the image of $\fie_D$ is a surface;
\item  if  $S$ has no irrational pencil, then  
 $\dim|D|\ge 2q-4$. 
  \end{enumerate}
\end{lem}
\begin{proof}  If $S$ has no irrational pencil of genus $q-1$, then the map $\fie_D$ is generically finite by Lemma \ref{gauss}.

Let $G^{*}\subset \pp(\bigwedge^2H^0(\Omega^1_S))$ be the Grassmannian of lines in $(\pp^{q-1})^*$. 
By the Theorem  of Castelnuovo--De Franchis, the surface $S$ has an irrational pencil if and only if there exist linearly independent  $\al, \be\in H^0(\Omega^1_S)$ such that  $v(\al\wedge \beta)=0$.
Assume that $S$ has no such pencil.  Then $\bar{v}$ restricts to a morphism $G^*\to \pp(\Lambda)$, which is finite onto its image, since $\bar{v}$ is induced by a linear map. It follows  $$\dim|D|=\dim\pp(\Lambda)\ge \dim G^*=2q-4.$$
\end{proof} 
\subsection{The cotangent map $\bar{\Phi}$.}

 The projectivized tangent bundle $\pp(\Omega^1_A)$ is canonically isomorphic to $A\times \pp^{q-1}$.  (Recall that, as explained in Notation and Conventions, for a vector space $V$ we denote by $\pp(V)$ the set  of $1$-dimensional subspaces of $V$, while for a vector bundle $\mathcal E$ we denote by $\pp({\mathcal E})$ the bundle $\Proj (\Sym {\mathcal E})$ of $1$-dimensional quotients).
Since the map $a$ is generically finite, its differential induces a rational map $\Phi\colon \pp(\Omega^1_S)\to A\times \pp^{q-1}$. Denoting by $\pi\colon \pp(\Omega^1_S)\to S$ and $p\colon A\times \pp^{q-1}$ the natural  projections, we have a commutative diagram:
\begin{equation}\label{diag1}
\begin{CD}
 \pp(\Omega^1_S)@>{\Phi}>>A\times \pp^{q-1} \\
@V{\pi}VV @VV{p}V\\
S @>{a}>> A 
\end{CD}
\end{equation} We denote by ${\bar{\Phi}}\colon \pp(\Omega^1_S)\to \pp^{q-1}$ the map obtained by  composing $\Phi$ with the projection $A\times \pp^{q-1}\to\pp^{q-1}$. The map $\bar{\Phi}$ is called the {\em cotangent map} (cf.  \cite{roulleau}). 
Clearly $\bar{\Phi}$ is the map given by  the linear system $|L|:=|\OO_{\pp(\Omega^1_S)}(1)|$. 
Let $E$ be the maximal effective divisor  of $S$ contained in the scheme of zeros of all the  $\eta\in H^0(\Omega^1_S)$. The fixed part of $|L|$ is $\pi^*E$.  We write $|L|=\pi^*E+|H|$ and for  $\eta \in H^0(\Omega^1_S)$ we denote the corresponding divisor by $L_{\eta}=\pi^*E+H_{\eta}$.
\begin{lem}\label{dim3}Assume that $q\ge 4$. 
If  the image $X$ of $\bar{\Phi}$ has dimension $<3$, then  $S$ has an irrational pencil  of genus $q-1$.
\end{lem}
\begin{proof}  By construction, $X$ is ruled by lines and the base of the ruling is birational to $\ga(\Si)$, where $\ga$ is the Gauss map. Since $q\ge 4$ by assumption and $X\subset \pp^{q-1}$ is non degenerate,  $X$ is not a plane. Since the plane is the only projective surface with a $2$-dimensional family of lines, it follows that $X$  has dimension $<3$  only if $\ga(\Si)$ is not a surface. By Lemma \ref{gauss} this happens  only if $S$ has an irrational pencil of genus $q-1$.
\end{proof}

\subsection{The map $\fie_{\eta}$.}

Let $\eta\in H^0(\Omega^1_S)$ and    denote by $\Lambda_{\eta}\subseteq H^0(K_S)$ the image of the map $H^0(\Omega^1_S)\to H^0(K_S)$ defined by $\al\mapsto \eta\wedge \al$ (cf. Proposition  \ref{xiaolines}). We write $\pp(\Lambda_{\eta})\!:=F_{\eta}+|D_{\eta}|$, where $|D_{\eta}|$ is the moving part. We denote by $\fie_{\eta}$ the map given by $|D_{\eta}$.

\begin{lem}\label{DEta} 
Let   $\eta\in H^0(\Omega^1_S)$ be general. Then:
\begin{enumerate}
\item $\dim |D_{\eta}|=q-2$;
\item $F_{\eta}=F$;
\item if $x$ is a base point (possibly infinitely near) of $|D|$  of multiplicity $m$, then the multiplicity of $x$ as a base point of  $|D_{\eta}|$ is also equal to $m$;
\item if $x$ is a base point of $|D_{\eta}|$ such that the differential of $a$ is non singular at $x$, then $x$ is a simple base point.
\end{enumerate} 
\end{lem}
\begin{proof}
By the Castelnuovo--De Franchis Theorem,  if $\al$ and $\eta$ are linearly independent 1-forms the form $\al\wedge \eta$ vanishes identically if and only if  there exists an irrational pencil $f\colon S\to B$ such that  $\al$ and $\eta$ are both pull backs of $1$-forms of $B$. Hence, for a general choice of $\eta\in H^0(\Omega^1_S)$ the system $|D_{\eta}|$ has dimension $q-2$.

Statement (ii) is \cite{manetti}, Lemma 1.2. Statement (iii) can be proven by the same argument. Namely,  given forms $\al,\be\in H^0(\Omega^1_S)$ with $\al\wedge\be\ne 0$, we write  the divisor of zeroes of $\al\wedge\be$ as $F+C_{\al,\be}$. Assume that  $x$ is an actual point of $S$: then $m$ is just the minimum of  the multiplicities of  the divisors $C_{\al,\be}$ at $x$, and the statement follows. If the point $x$ is infinitely near, then one applies the same argument on a suitable blow up of $S$.

Let $U\subseteq S$ be the open set of points where the differential of $a$ is non singular. A point $x\in U$  is a base point of $|D_{\eta}|$ iff and only if $\eta(x)=0$. The bundle $\Omega^1_S|_U$ is generated by the elements of $H^0(\Omega^1_S)$, hence the zeros of a general $\eta$ on $U$ are simple. It is easy to check that this implies that they are also simple base points of $|D_{\eta}|$.
\end{proof} 

\begin{lem}\label{fieEta}

Assume that     $\eta\in H^0(\Omega^1_S)$ is  general.
\begin{enumerate}
\item If $q\ge 4$ and $S$ has no irrational pencil of genus $q-1$, then  the image of $\fie_{\eta}$ is a surface;
\item if $q\ge 5$ and  $S$ has no irrational pencil, then the image of $\fie_{\eta}$ is not ruled by lines.
\end{enumerate}
\end{lem}
\begin{proof}
To a general  $1$-form $\eta\in H^0(\Omega^1_S)$ there corresponds a general  element $L_{\eta}$ of $|L|:=|\OO_{\pp(\Omega^1_S)}(1)|$. As we have already observed, we can write $L_{\eta}=\pi^*E+H_{\eta}$, where $E$ is the divisorial part of the base locus of all the global $1$-forms of $S$ and $H_{\eta}$  is irreducible. Thus $H_{\eta}$ defines a rational section $\si_{\eta}\colon S\to \pp(\Omega^1_S)$ and $\fie_{\eta}$ is just the composite map $\bar{\Phi}\circ \si_{\eta}$, where $\bar{\Phi}$ is the cotangent map. So the image of $\fie_{\eta}$ is a general hyperplane section of the image $X$ of $\bar{\Phi}$ and statement (i)  follows by Lemma \ref{dim3}.

To prove  (ii), assume for contradiction that  the image of $\fie_{\eta}$ is ruled by lines. 
Since $S$ has no irrational pencil, the pull back $|F|$ of this ruling  is either a linear pencil or an elliptic pencil. In either case Proposition \ref{xiaolines} gives $$q-3\le\dim|D_{\eta}|_F= 1,$$ a contradiction. \end{proof}

\begin{cor}\label{Dlines}
%Let $S$ be a surface of maximal Albanese dimension with irregularity 
If $q\ge 5$ and  $S$ has no irrational pencil, then the image of the canonical map $\fie$ and the image  of $\fie_D$ are  not ruled by lines.
\end{cor}
\begin{proof}
The claim follows immediately from Lemma \ref{fieEta}, since $\fie_{\eta}$ is obtained by  composing $\fie$ or  $\fie_D$ with a projection. \end{proof}

\subsection{Some inequalities}
\begin{defn}  Let $|M|$ be a linear system without fixed divisors on a surface $S$,  let $\fie_M$ be the rational map given by $|M|$ and let $Y$ be the image of $\fie_M$. We define the {\em free self intersection number of $|M|$}, denoted by $\fdeg|M|$, as follows:
\begin{itemize}
\item $\fdeg|M|=deg\fie_M\deg Y$ if $Y$ is a surface;
\item $\fdeg|M|=0$ otherwise.
\end{itemize}
\end{defn}
\smallskip

\begin{remark}  Let  $\epsi\colon S'\to S$ be a blow up that solves the indeterminacy of $\fie_M$ and let $|M'|$ be the strict transform of $|M|$ on $S'$. Then $\fdeg|M|=(M')^2$.
\end{remark}

\begin{prop}\label{fdegD}
%Let $S$ be a surface of maximal Albanese dimension, with irregularity $q\ge 5$.
Assume that  $S$ has no irrational pencil and $q\ge 5$.
\begin{enumerate}
\item If $\fie_D$ is birational then, $\fdeg|D|\ge 6q-15$;
\item if $\fie_D$ is not birational, then $\fdeg|D|\ge 8q-25$.
\end{enumerate}
\end{prop}
\begin{proof}
By the  remark above, up to replacing $S$ by a suitable blow up, we may assume that the system $|D|$ is free, so  that  $\fie_D$ is a morphism and $D^2=\fdeg|D|$.

Set $r=\dim|D|$ and denote by $\Si_D\subset \pp^r$ the image of $\fie_D$. Recall that by Lemma \ref{fieD} $\Si_D$ is a surface and $r\ge 2q-4$. 

Assume that $\fie_D$ is birational. A general $D$ is smooth and  the system $|D|_D$ is special on   $D$. By Lemma \ref{clifford} we have:
\[
\fdeg |D|=D^2\ge 3h^0(\OO_D(D))-4\ge 3r-4\ge 6q-16
\]
Assume now that $\fdeg |D|=6q-16$ and write $K_S=F+D$, with $F\ge 0$.  By Lemma \ref{clifford}, we have $\OO_D(2D)=K_D$, hence $FD=0$.  Let $\bar S$ be the minimal model of  $S$,  let $\bar D$ be the image of $D$ and $\bar F$ the image of $F$. Since $F$ contains the fixed part of $|K_S|$ and $DF=0$, the system $|\bar D|$ is free and  we have 
 $\bar F\bar D=0$, hence $\bar F=0$ by the $2$-connectedness of canonical divisors (the assumptions imply that $S$ is of general type). So we have $\fdeg |D|=\bar D^2=K_{\bar S}^2$.  By \cite[Thm.3.2]{deb1} we get:
\[
\fdeg|D|=K^2_{\bar S}\ge3 p_g(S)-7+q \ge 3r-4+q> 6q-16,
\]
 a contradiction.   This  completes the proof of statement (i).

Assume that $\fie_D$ is not birational. If $\deg\fie_D=2$, then $p_g(\Si_D)>0$ by Lemma \ref{involution}. Hence $\Si_D$ has degree $\ge 2r-2$ by Proposition \ref{ruled}. So we have:
$$\fdeg|D|\ge 2(2r-2)\ge 8q-20.$$
If $\deg\fie_D\ge 3$ and $q\ge 5$, then by Corollary \ref{Dlines} the surface $\Si_D$ is not ruled by lines. By Proposition \ref{Mileslines}, we have:
$$\fdeg|D|=D^2\ge 4(2q-4)-9=8q-25.$$
\end{proof}

Next we give analogous inequalities for the map $\fie_{\eta}$.
\begin{prop}\label{fdegEta}
%Let $S$ be a surface of maximal Albanese dimension, with irregularity $q\geq 5$
If $q\ge 5$,  $S$ has no irrational pencil and $\eta\in H^0(\Omega^1_S)$ is  general, then:
\begin{enumerate}
\item if $\fie_{\eta}$ is birational then, $\fdeg|D_{\eta}|\ge 3q-9$;
\item if $\fie_{\eta}$ is not birational, then $\fdeg|D_{\eta}|\ge 4q-17$.
\end{enumerate}
\end{prop}
\begin{proof} The proof follows verbatim the proof of Proposition \ref{fdegD}, recalling that here $r=q-2$ by Lemma \ref{DEta}.
\end{proof}

Finally, we give an estimate for the geometric genus of a general $D$.
\begin{prop}\label{gD}
%Let $S$ be a surface of maximal Albanese dimension, with irregularity  $q\geq 5$
If $q\ge 5$,   $S$ has no irrational pencil and $D$ is a general curve of $|D|$. Then:
\begin{enumerate}
\item if $\fie_D$ is not birational, then $g(D)\ge 8q-24$;
\item if $\fie_D$ is birational, then $g(D)\ge 7q-15$.
\end{enumerate}
\end{prop}
\begin{proof}
Up to blowing up $S$, we may assume that $D$ is smooth and  $D^2=\fdeg|D|$. Then statement (i) follows directly from Proposition \ref{fdegD} by using adjunction on $S$.

If $\fie_D$ is birational,  then $h^0(2D)\geq 4h^0(D)-6$ by \cite{deb1}, Prop. 3.1.  So the image of  the restriction map  $H^0(2D)\to H^0(\OO_D(2D))$ has dimension  $\geq 3h^0(D)-6$. This implies that the image of the restriction map $r:H^0(K_S+D)\to H^0(K_D)$ has dimension  $\geq 3h^0(D)-6\ge 6q-15$.  The divisor $D$ is nef and big (cf. Lemma \ref{fieD}), 
hence  $H^1(K_S+D)=0$ by Kawamata--Viehweg vanishing. It follows $g(D)-q\ge 6q-15$.\end{proof} 

\begin{cor}\label{unramif} 
Assume that $q\ge 6$,  $S$ has no irrational pencil and  $\eta \in H^0(\Omega^1_S)$ is  general.  
If $a$ is unramified in codimension 1 and $\fie_{\eta}$ is birational, then $$\fdeg |D_{\eta}|\ge 4q-15.$$
\end{cor}
\begin{proof}
By Proposition \ref{DEta}, the geometric genus of $D_{\eta}$ is equal to the geometric genus of $D$, hence $\ge 7q-15$ by Proposition \ref{gD}. 

 Let $d=\fdeg|D_{\eta}|$, and write $d-1=m(q-4)+\epsi$, where $m\ge 0$ and $0\le \epsi\le q-5$ are integers.  The image via $\fie_{\eta}$ of a general $D_{\eta}$ is a curve of degree $d$ in $\pp^{q-3}$. By Castelnuovo's  bound we have:
  $$7q-15\le g(D)\le {m \choose 2} (q-4)+m\epsi.$$
  This implies $m\ge 4$, hence $d\ge 4q-15$.
\end{proof}

\section{Proofs of Theorem \ref{main}, \ref{main3} and \ref{main2}.}\label{proofs}

\begin{proof}[Proof of Theorem \ref{main3}]
Let $\fie$ be the canonical map of $S$ and let $\Si$ be the canonical image.
If $\fie$ has degree 2 and $S$ has no irrational pencil, then by Lemma \ref{involution}, $p_g(\Si)>0$. Then by \cite{beauville}, Thm. 3.1, the surface $\Si$ is canonically embedded in $\pp^{p_g(S)-1}$, hence by \cite{beauville}, Thm. 5.5, it has degree $\ge 3p_g(S)-7$. 
So we have:
$$K^2_S\ge 6p_g(S)-14=6\chi(S)+6q-20.$$

If $\fie$ has degree $\ge 3$, then by Proposition \ref{mainpencil} it is enough to consider the case when $S$ has no irrational pencil. In this case   $\Si$ is  not ruled by lines  by  Corollary \ref{Dlines}. Proposition \ref{Mileslines} then gives:
$$K^2_S\ge 4p_g(S)-13=4\chi(S)+4q-17.$$
\end{proof}
\begin{proof}[Proof of Theorem \ref{main}]
If $S$ has an irrational pencil, then the claim follows by the  stronger inequalities of Proposition  \ref{mainpencil}.
 Also, we may assume that $K_S-D>0$, since otherwise we have the stronger inequality of Theorem \ref{main2}, (ii).

Otherwise we prove the claim by refining  the estimate of \cite[Proof  of 6.2]{manetti}. We  recall  the notation of \cite{manetti} and the outline of the proof. We write $|L|:=|\OO_{\pp(\Omega^1_S)}(1)|$. Let $E$ be the maximal effective divisor  of $S$ contained in the scheme of zeros of every $\eta\in H^0(\Omega^1_S)$. The fixed part of $|L|$ is $\pi^*E$, where $\pi\colon \pp(\Omega^1_S)\to S$ is the projection.  We write $|L|=\pi^*E+|H|$. One observes that 
 \begin{gather*}
 3(K^2_S-4\chi(S))=L^2(L+\pi^*K_S)=\\
 =LH(L+\pi^*K_S)+2K_SE\ge LH(L+\pi^*K_S),
 \end{gather*}
 where the inequality is a consequence of the fact that $S$ is minimal.

The proof consists in estimating the right hand term of the above inequality by means of  geometric arguments. 
Set, as usual, $F:=K_S-D$ with $F>0$  and write $F=E+\sum \be_jB_j+\ga_lG_l$, where $\be_j$, $\ga_l$ are positive integers,  the $B_j$ are the curves contracted by $a$ and the $G_j$ are the  divisorial components of the ramification locus of $a$ that are not contracted.
For  $\eta \in H^0(\Omega^1_S)$ we denote the corresponding divisor of $|L|$ by $L_{\eta}=\pi^*E+H_{\eta}$.
 Let $\eta,\mu\in H^0(\Omega^1_S)$ be general. The divisor of zeros of $\eta\wedge \mu$ can be written as $F+D_{\eta,\mu}$, where $D_{\eta,\mu}$ is general in $|D_{\eta}|$. Consider the effective $1$-cycle $[L_{\eta}\cap H_{\mu}]$. One has:
 \[
 [L_{\eta}\cap H_{\mu}]=\sum_i m_iF_{x_i}+\sum_j \be_jR_j+\sum_l\ga_lS_l+C,
 \]
 where the $m_i$ are positive integers, $F_{x_i}$ is the fibre of $\pi$ over the point $x_i$,  the projection $\pi$ maps  birationally $R_j$ onto $B_j$, $S_l$ onto $G_l$ and $C$ onto $D_{\eta,\mu}$.
We set:
\begin{gather*}
A_0=(L+\pi^*K_S)(\sum_i m_iF_{x_i}+\sum_j \be_jR_j),\\
A_1=(L+\pi^*K_S)\sum_l\ga_lS_l,
\\A_2=(L+\pi^*K_S)C,
\end{gather*}
 so that: $$3(K^2_S-4\chi(S))\ge(L+\pi^*K_S)[L_{\eta}\cup H_{\mu}]=  A_0+A_1+A_2.$$
 Then we have:
 \begin{equation*}
 A_2=LC+K_SD_{\eta,\mu}\ge HC+K_SD_{\eta,\mu}= HC+D^2+D(K_S-D).
 \end{equation*}
 In \cite[Prop. 3.5]{manetti} it is shown that for every $l$ one has $$LS_l\ge -(G_l)^2.$$  This gives:
 \begin{equation*}
 A_1\ge \sum_l\ga_l(K_S-G_l)G_l\ge (K_S-\ga_lG_l)G_l,
 \end{equation*}
  where the last inequality follows from the fact that $K_S$ is nef.
In addition, Manetti shows that:
 \begin{equation*}
 A_0\ge-\si,
 \end{equation*}
where $\si$ is the number of divisors $0\le Z\le F$ such that:

1)  the support of $Z$ is a connected component of $\cup_jB_j$;

2) all the components of $Z$ are $-3$-curves;

3) $Z$ is maximal among divisors with properties 1) and 2).

Assume that there is such a divisor $Z$ with $Z>0$.   Notice that, since the support of $Z$ is contracted by $a$,  the intersection form on the components of $Z$ is negative definite, so that, in particular, $Z^2<0$. The number  $Z(K_S-Z)$ is  $\ge 2$ by the $2$-connectedness of canonical divisors and it is even by the adjunction formula.  If $Z=\Ga$, where $\Ga$ is an  irreducible $-3$-curve, then $Z(K_S-Z)=4$. If $Z$ is not reduced or has more than one component, then $Z(K_S-Z)\ge 2-Z^2\ge 3$, hence  $Z(K_S-Z)\ge 4$. It follows that:
\[
D(K_S-D)+\sum_l(K_S-\ga_lG_l)G_l-\si\ge D(K_S-D)/2+1.
\]
The inequality is also trivially satisfied if there is no $Z$ with properties 1),2), and 3) (i.e.,   if $\si=0$), since we have assumed that $K_S-D>0$.
Thus:
\[
D^2+D(K_S-D)+\sum_l(K_S-\ga_lG_l)G_l+A_0\ge p_a(D).
\]
All the previous estimates together give:
\[ 
3(K^2_S-4\chi(S))\ge HC+p_a(D).
\]
Since $HC\ge\fdeg|D_{\eta}|$,   by Propositions \ref{fdegEta} and \ref{gD}  we have
\[ K^2_S\ge 4\chi(S)+\frac{10}{3}q-8.
\]

\end{proof}
\begin{proof}[Proof of Theorem \ref{main2}]
Consider statement (i) first. One argues exactly as in the proof of Theorem \ref{main}. The only difference is that,  since $a$ is not birational, both $\fie_D$  and $\fie_{\eta}$ are not birational either and Propositions \ref{fdegEta} and \ref{gD} give stronger inequalities in this case.

Consider now statement (ii). The assumption on the Albanese map is equivalent to the fact that $L$ is nef on $\pp(\Omega^1_S)$ or, equivalently, $K_S-D=0$. Arguing as in the proof of Theorem \ref{main} and using the same notation, we have:
\[ 
2(K^2_S-6\chi(S))=L^3\ge LC\ge \fdeg |D_{\eta}|.
\]
For $q\ge 6$, the statement now follows by Proposition \ref{fdegEta} and Corollary \ref{unramif}.
For $q=5$, and $\fie_{\eta}$ not birational we use Proposition \ref{fdegEta} again, while if $\fie_{\eta}$ is birational we have $\fdeg|D_{\eta}|\ge 4$ since the curve $D_{\eta}$ is of general type.

\end{proof}

\bigskip

\bigskip

\begin{minipage}{13cm}
\parbox[t]{6.5cm}{Margarida Mendes Lopes\\
Departamento de  Matem\'atica\\
Instituto Superior T\'ecnico\\
Universidade T{\'e}cnica de Lisboa\\
Av.~Rovisco Pais\\
1049-001 Lisboa, PORTUGAL\\
mmlopes@math.ist.utl.pt
 } \hfill
\parbox[t]{5.5cm}{Rita Pardini\\
Dipartimento di Matematica\\
Universit\`a di Pisa\\
Largo B. Pontecorvo, 5\\
56127 Pisa, Italy\\
pardini@dm.unipi.it}
\end{minipage}


\begin{thebibliography}{ABCD}

\bibitem[BNP]{bnp} M.A.~Barja, J.C,~Naranjo, G.P.~Pirola, {\em On the topological index of irregular surfaces},  J. Algebraic Geom.  {\bf 16}  (2007),  no. 3, 435--458.
\bibitem[Be1]{beauville} A. Beauville, {\em L'application canonique pour les surfaces de type g\'en\'eral}. Inv. Math. {\bf 55} (1979), 121--140.
\bibitem[Be2]{appendix} A. Beauville, {\em L'in\'egalit\'e $p_g\ge 2q-4$ pour les surfaces de type g\'en\'eral}, appendix to \cite{deb1}.

\bibitem[Ca]{ca}  F. Catanese, {\em On the moduli spaces of surfaces of general type}, J. Diff. Geometry {\bf 19} (1984), 483--515.
\bibitem[CH]{ch} M. Cornalba, J. Harris, {\em Divisor classes associated to families of stable varieties, with applications to the moduli space of curves}, Ann. Sci. \'Ecole Norm. Sup., (4) {\bf 21} (1988), no. 3, 455--475.

\bibitem[De1]{deb1} O. Debarre, {\em In\'egalit\'es num\'eriques pour les surfaces de type g\'en\'eral}, with an appendix by A. Beauville, Bull. Soc. Math. France {\bf 110} 3 (1982),  319--346.
\bibitem[De2]{deb2} O. Debarre, {Tores et vari\'et\'es ab\'eliennes}, S.M.F. Sciences 1999.
\bibitem[Ha]{harris} J. Harris, {\em A bound on the geometric genus of projective varieties}, 
Ann. Scuola Norm. Sup. Pisa Cl. Sci. (4) 8 (1981), no. 1, 35--68.

\bibitem[Ko]{konnoeven} K. Konno, {\em Even canonical surfaces with small $K^2$} III,  Nagoya Math. J.  {\bf 143} (1996), 1--11.
\bibitem[Ma]{manetti} M. Manetti, {\em Surfaces of Albanese general type and the Severi conjecture}, Math. Nach. {\bf 261--262} (2003), 105-122.
\bibitem[Pa]{severi} R. Pardini, {\em The Severi inequality $K^2\ge 4\chi$ for surfaces of maximal Albanese dimension},  Invent. math. {\bf 159}  3 (2005), 669 --672. 
\bibitem[Re1]{milesdeg} M. Reid, {\em Surfaces of small degree}, Math. Ann. (1986), 71--80.
\bibitem[Re2]{laquila} M. Reid, {\em Quadrics through a canonical surface}, in: Algebraic Geometry -- Hyperplane
sections and related topics (L'Aquila 1988), Springer LNM 1417 (1990), 191--213.
%\bibitem[Xi2]{xiaohyp} G. Xiao, {\em Hyperelliptic surfaces of general type with $K^2<4\chi$}, Manuscripta Math. {\bf 57} (1987), 125--148.
\bibitem[Rou]{roulleau} X. Roulleau, {\em L'application cotangente des surfaces de type g\'en\'eral}, Th\`ese de doctorat, Universit\'e d'Angers (2007).
\bibitem[St]{stoppino} L. Stoppino, {\em Slope inequalities for fibred surfaces via GIT}, Osaka Journal of Mathematics, Vol. 45, No. 4 (2008)
\bibitem[Xi1]{xiaopencil} G.~Xiao {\em Irregularity of surfaces with a linear pencil}  Duke Math. J.  {\bf 55}  (1987),  no. 3, 597--602.
\bibitem[Xi2]{xiaoslope} G. Xiao, {\it Fibered algebraic surfaces with low slope}, Math. Ann. {\bf 276}
(1987), 449--466.
\end{thebibliography}
\end{document}